\documentclass[leqno]{article}
\usepackage[active]{srcltx}
\usepackage{amssymb,amsmath}
\newtheorem{proposition}{Proposition}[section]
\newtheorem{theorem}{Theorem}[section]

\newtheorem{corollary}{Corollary}[section]
\newtheorem{definition}{Definition}[section]
\newtheorem{remark}{Remark}[section]

\numberwithin{equation}{section}

\begin{document}

\begin{center}

{\Large A Nash type solution for hemivariational inequality systems
Systems}

\vspace*{1.5cm}

 {\bf Du\v san Repov\v s  and  Csaba Varga}
\medskip

 {\small $^c\,$Faculty of Mathematics and Physics,
and Faculty of Education, University of Ljubljana,\\
\small  P. O. B. 2964, Ljubljana, Slovenia 1001\\
Babe\c s-Bolyai University, Faculty of Mathematics and Computer
Sciences, RO-400084 Cluj-Napoca, Romania}\\ {\it Email address}:
{dusan.repovs@guest.arnes.si \ and \ csvarga@cs.ubbcluj.ro}
\end{center}


\vspace{1cm}
\begin{abstract}
 In this paper we prove an existence result for a general
class of hemivariational inequalities systems using the Ky Fan
version of the KKM theorem (1984) or the Tarafdar fixed point theorem
(1987). As application we give an  infinite dimensional
version for existence result of Nash generalized derivative points
introduced recently by Krist\'{a}ly (2010) and also we
give an application to a general hemivariational inequalities
system.
\end{abstract}

\bigskip



 \section{Introduction}

In the last years many papers have been dedicated to the of study
the existence and multiplicity of solutions for hemivariational
inequality systems or differential inclusion systems defined on
bounded or unbounded domain, see \cite{BV},\cite{BHV},
\cite{Kristaly1}, \cite{Kristaly2}, \cite{Kristaly3}, \cite{LV}.
In these papers the authors use critical points theory for locally
Lipschitz functions, combined with the {\it Principle of Symmetric
Criticality} and different topological methods.  For a
comprehensive treatment of hemivariational inequality and
hemivariational inequalities systems on bounded domains using
critical point theory for nonsmooth functionals we refer to the
monographs of  D. Motreanu and V. R\u adulescu \cite{MR} and D.
Motreanu and P. D. Panagiotopoulos \cite{MP}. For very recent
results concerning variational inequalities and elliptic systems
using critical points theory and different variational methods see
also the book of Krist\'{a}ly, R\u{a}dulescu and Varga
\cite{KRV}.\\

 The aim of this paper is to prove that the
 existence of at least one solution for a general class of
 hemivariational inequalities systems on  a
 closed and convex set (either bounded or unbounded)
 without using critical point theory.
 We apply a version of the well-known theorem of Knaster--Kuratowski--Mazurkiweicz
 due to Ky Fan \cite{Fan_fix} or the Tarafdar fixed point theorem
 \cite{TarFix}.
 We start the paper by giving in Section 2 the assumptions  and formulating the
 hemivariational inequalities system problem which we study.
  The main results concerning the existence of at least one solution for the hemivariational
 inequalities systems which we study are given in Section 3. Section 4 contains applications to Nash
 and Nash general derivative points  and existence results for some abstract class of hemivariational inequalities
 systems.

 \section{Assumptions and Formulation of the Problem}

Let $X_{1}, X_{2}, \ldots, X_{n}$ be reflexive Banach spaces and
$Y_{1}, Y_{2}, \ldots, Y_{n}, Z_{1}, \ldots, Z_{n}$ Banach spaces
such that there exist linear operators $T_{i}: X_{i} \rightarrow
Y_{i}$, $T_{i}: X_{i} \rightarrow Z_{i}$ for $i\in \{1, \ldots,
n\}$. We suppose that the following condition
hold:\\
{\bf (TS)} $T_{i}: X_{i} \rightarrow Y_{i}$ and $S_{i}:X_{i}
\rightarrow Z_{i}$ are compact for $i=\overline{1,n}$.

We denote by $X_{i}^{*}$ the topological dual of $X_{i}$ and
$\langle\cdot,\cdot\rangle_{i}$ denotes the duality pairing
between $X_{i}^*$, whereas  $X_{i}$ for $i=\overline{1, n}$.
 Also, let $K_{i} \subset X_{i}$ be closed, convex sets for
$i=\overline{1,n}$ and  we consider $A_{i}: Y_{1} \times \cdots
\times Y_{i}\times \cdots \times Y_{n} \rightarrow \mathbb{R}$ the
continuous functions which are locally Lipschitz in the $i$th
variable and we denote by $A_{i}^{\circ}(u_{1}, \ldots, u_{i},
\ldots, u_{n}; v_{i})$  the partial Clarke derivative in the
directional derivative in the $i$th variable, i.e. the Clarke
derivative of the locally Lipschitz function $A_{i}(u_{1}, \ldots,
u_{i}, \ldots, u_{n})$ at the point $u_{i}\in Y_{i}$ in the
direction $v_{i} \in Y_{i}$, that is
$$A_{i}^{\circ}(u_{1}, \ldots, u_{i}, \ldots u_{n};v_{i})=
\limsup_{\scriptstyle {\it w\to u_{i}}\atop
 \scriptstyle \it \tau\searrow 0}\frac{A_{i}(u_{1},
\ldots,w+ \tau v_{i}, \ldots u_{n})- A_{i}(u_{1}, \ldots,w, \ldots
u_{n})}{\tau}
$$
We suppose that for every $i=\overline{1,n}$ the following
condition holds:\\
{\bf (A)} the functions $A_{i}^{\circ}: Y_{1} \times \cdots \times
Y_{n} \times Y_{i}\rightarrow \mathbb{R}$ are upper
semi-continuous.

We also consider the following nonlinear operators $F_{i}:K_{1}
\times \cdots \times K_{i} \times \cdots K_{n} \rightarrow
X_{i}^{\star}, i=\overline{1,n}$. We suppose that the operators
$F_{i}$ satisfy the
following condition:\\
{\bf(F)} the functions $(u_{1}, \ldots, u_{n}) \mapsto \langle
F_{i}(u_{1}, \ldots, u_{n}), v_{i}\rangle_{i}$ are weakly upper
semi-continuous for every $v_{i} \in X_{i}$ and
$i=\overline{1,n}$.

\begin{definition}\label{regular_lip}{\rm [see \cite{Cl}]}
Let $Z$ be a Banach space and  $j: Z \rightarrow \mathbb{R}$  a
locally Lipschitz function. We say that $j$ is regular at $u \in
Z$ if for all $v \in Z$ the one-sided directional derivative
$j^{\prime}(u;v)$ exists and $j^{\prime}(u;v)=j^{\circ}(u;v)$. If
$j$ is regular at every point $u \in Z$ we say that $j$ is
regular.
\end{definition}

We have the
following elementary result.

\begin{proposition}\label{Regular_prop}
Let  $J: Z_{1} \times \cdots \times Z_{n} \rightarrow \mathbb{R}$
be a regular, locally Lipschitz function. Then the following
assertions hold:
\begin{itemize}
\item [a)] $\partial J(u_{1}, \ldots, u_{n}) \subseteq \partial_{1}J(u_{1}, \ldots, u_{n}) \times \cdots \times \partial_{n} J(u_{1}, \ldots, u_{n})$
(see {\rm\cite{Cl}, Proposition 2.3.15)}, where $\partial_{i} J,
i=\overline{1,n}$ denotes the Clarke subdifferential in the $i$th
variable;
\item [b)] $J^{\circ}(u_{1}, \ldots, u_{n};v_{1}, \ldots, v_{n} ) \leq \displaystyle \sum_{i=1}^{n} J_{i}^{\circ}(u_{1}, \ldots,
u_{n};v_{i})$, where $J_{i}^{\circ}$ denotes the Clarke derivative
in the {\rm i}th variable; and
\item [c)] $J^{\circ}(u_{1}, \ldots, u_{n}; 0, \ldots, v_{i}, \ldots,
0) \leq J_{i}^{\circ}(u_{1}, \ldots, u_{n}; v_{i}) $.
\end{itemize}
\end{proposition}

We introduce the following notations:
\begin{itemize}
\item   $K=K_{1}
\times \cdots \times K_{n},$

\item   $u=(u_{1}, \ldots, u_{n})$

\item   $Tu=(T_{1}u_{1}, \ldots, T_{n}u_{n})$
\item   $Su=(S_{1}u_{1}, \ldots, S_{n}u_{n})$
\item   $A(Tu, Tv-Tu)= \displaystyle \sum_{i=1}^{n}A_{i}^{\circ}(Tu, T_{i}v_{i}- T_{i}u_{i} )$

\item   $F(u, v-u)=\displaystyle \sum_{i=1}^{n} \langle F_{i}u,
v_{i}-u_{i}\rangle_{i}$.
\end{itemize}
In this paper we study the following problem:\\
{\it Find $u=(u_{1}, \ldots, u_{n}) \in K_{1}\times\ldots \times
K_{n}$ such that for all $v=(v_{1}, \ldots, v_{n}) \in
K_{1}\times\ldots \times K_{n}$ and $i \in \{1, \ldots, n\}$ we
have:
$$A_{i}^{\circ}(Tu;Tv_{i}-Tu_{i}) +
\langle F_{i}(u), v_{i}-u_{i} \rangle_{i} + J^{\circ}_{i}(Su;
S_{i}v_{i}- S_{i}u_{i})\geq 0. \leqno{{\bf (QHS)}}$$ In this case
we say that $u=(u_{1}, \ldots, u_{n})$ is a Nash equilibrum point
for the system {\bf (QHS)}}.

To prove our main result we use the FKKM theorem due to Ky Fan
\cite{Fan_fix} and the Tarafdar fixed point theorem
\cite{TarFix}.

\begin{definition}\label{kkm_ertelmezes}
Suppose that $X$ is a vector space and $E\subset X$. A set-valued
mapping $G:E\to 2^X$ is called a KKM mapping, if for any
$x_1,\dots, x_n \in E$ the following holds
$$ {\rm conv}\{x_1,\dots, x_n \} \subset \bigcup_{i=1}^n G(x_i).$$
\end{definition}

The following version of the KKM theorem is due to Ky Fan
\cite{Fan_fix}.

\begin{theorem} \label{kkm}
Suppose that $X$ is a locally convex Hausdorff  space, $E\subset
X$ and that $G:E\to 2^X$ is a closed-valued KKM map. If there
exists $x_0\in E$ such that $G(x_0)$ is compact, then
$\displaystyle \bigcap_{x\in E} G(x)\neq \emptyset.$
\end{theorem}

\begin{theorem}\label{Tarafdar_Fix}
Let $K$ be a nonempty, convex subset of a Hausdorff topological
vector space X. Let $G: K \hookrightarrow 2^{K}$ be a setvalued
map such that
\begin{itemize}
\item [i)] for each $u \in K$, $G(u)$ is a nonemty convex subset
of $K$;
\item [ii)] for each $v \in K$, $G^{-1}(v)=\{ u \in K \ \ : \ v \in
G(u)\}$ contains an open set $O_{v}$ which may be empty;
\item [iii)] $\displaystyle \cup_{v \in K}O_{v} =K$; and
\item [iv)] there exists a nonemty set $K_{0}$ contained in a
compact convex subset $K_{1}$ of $K$ such that $D=\displaystyle
\cap_{v \in K_{0}} O_{v}^{c}$ is either empty or compact (where
$O_{v}^{c}$ is the complement of $O_{v}$ in $K$).
\end{itemize}
Then there exists a point $u_{0} \in K$ such that $u_{0} \in
G(u_{0})$.

\end{theorem}

E. Tarafdar in \cite{TarFix} proved the equivalence of Theorems
\ref{kkm} and \ref{Tarafdar_Fix}.

\section{Main results}

\begin{theorem}\label{Main_RES} Let $K_{i} \subset X_{i},
i=\overline{1,n}$ be nonempty, bounded, closed and convex sets.
Let $A_{i}: Y_{1} \times \cdots \times Y_{i}\times \cdots \times
Y_{n} \rightarrow \mathbb{R}$ be a locally Lipschitz function in
the $i$th variable for all $i \in \{1, \ldots, n\}$ satisfying
condition {\bf(A)}. We suppose that the operators $T_{i}:
X_{i}\rightarrow Y_{i}$, $S_{i}: X_{i}\rightarrow Z_{i}$ and
$F_{i}: K_{1} \times \ldots \times K_{n} \rightarrow
X_{i}^{\star}$($i=\overline{1,n}$) satisfy the condition {\bf(TS)}
respectively {\bf (F)}. Final we consider the regular locally
Lipschitz function $J: Z_{1} \times \cdots \times Z_{n}
\rightarrow \mathbb{R}$.  Under these conditions the problem (QHS)
admits at least one solution.
\end{theorem}

Before proving Theorem \ref{Main_RES}, we make two remarks.

\begin{remark}\label{remark1}
We observe that for every $v \in K$ the function
$$u \mapsto A(Tu, Tv-Tu) + F(u, v-u)+ J^{\circ}(Su;Sv-Su)$$ is weakly upper
semi-continuous. Indeed, from the condition {\bf (A)} and from the
fact that the operators $T_{i}$ are compact follows that $A(Tu,
Tv-Tu)$ is weakly upper semi-continuous. From {\bf (F)} it follows
that $F(u, v-u)$  is weakly upper semi-continuous. The third term,
i.e. $ J^{\circ}(Su;Sv-Su)$ is weakly upper semi-continuous,
because $J^{\circ}(\cdot;\cdot)$ is upper semi-continuous and the
operators $S_{i}: X_{i} \rightarrow Z_{i}$ are compact.
\end{remark}

\begin{remark}\label{remark2}
 If there exists $u \in K$,
such that for every $v \in K$ we have:
\begin{equation}\label{FO_egyenlet}
A(Tu, Tv-Tu) + F(u, v-u)+ J^{\circ}(Su;Sv-Su) \geq 0,
\end{equation}
then $u \in K$ is a solution of the problem (QHS).
Indeed, if we fix an $i=\{1, \ldots, n\}$ and put $v_{j}:=u_{j}, j
\neq i$ in the above inequality and using  iii) Proposition
\ref{Regular_prop} we get that
$$A_{i}^{\circ}(Tu;Tv_{i}-Tu_{i}) +
\langle F_{i}(u_{i}), v_{i}-u_{i} \rangle_{i} + J^{\circ}_{i}(Su;
S_{i}v_{i}- S_{i}v_{i})\geq 0. \leqno{(QHS)},$$ for all $i \in
\{1, \ldots, n\}$.
\end{remark}
In the sequel we give two proofs, using Theorems \ref{kkm} and
\ref{Tarafdar_Fix}.\\

{\bf First Proof:} Let $G: K \hookrightarrow 2^{K}$ be the
set-valued map defined by
$$G(v)= \{ u \in K \ : A(Tu, Tv-Tu) + F(u, v-u)+ J^{\circ}(Su;Sv-Su) \geq 0\}.$$
For every $v \in K$, we have $G(v) \neq \emptyset$ because $v \in
G(v)$ and taking into account that the function
$$u \mapsto A(Tu, Tv-Tu) + F(u, v-u)+ J^{\circ}(Su;Sv-Su)$$ is weakly upper
semi-continuous, it follows that the set $G(v)$ is weakly closed.
Now we prove that $G$ is a KKM mapping. We argue by contradiction,
let $v_{1}, \ldots, v_{k} \in K$ and $w \in {\rm conv} \{v_{1},
\ldots, v_{k}\}$ such that $w \notin \displaystyle \cup_{i=1}^{k}
G(v_{i})$. From this it follows that
\begin{equation}\label{1_egyenlet}
A(Tw, Tv_{i}-Tw) + F(w, v_{i}-w)+ J^{\circ}(Sw;Sv_{i}-Sw) < 0,
\end{equation}
 for
all $i=\{1, \ldots, k\}$. Because of $w \in {\rm conv}\{v_{1},
\ldots, v_{k}\}$ the existence of $\lambda_{1}, \ldots,
\lambda_{k} \in [0,1]$ with $\displaystyle \sum_{i=1}^{k}
\lambda_{i}=1$ such that $w=\displaystyle \sum_{i=1}^{k}
\lambda_{i}v_{i}$ follows . If we multiply the inequalities
(\ref{1_egyenlet}) with $\lambda_{i}$ and adding for $i=\{1,
\ldots, k\}$ we obtain

\begin{equation}\label{2_egyenlet}
A(Tw, Tw-Tw) + F(w, w-w)+ J^{\circ}(Sw;Sw-Sw) < 0
\end{equation}
because the functions $A(\cdot, \cdot), F(\cdot, \cdot)$ and
$J^{\circ}(\cdot, \cdot)$ are positive homogeneous and convex in
the second variable. From inequality (\ref{2_egyenlet}) it follows
that $0= A(Tw, Tw-Tw) + F(w, w-w)+ J^{\circ}(Sw;Sw-Sw) < 0$, which
is a contradiction. Because the set $K$ is bounded, convex and
closed, it follows that it is weakly closed and by the
Eberlein-Smulian theorem we have is weakly compact. Because $G(v)
\subset K$ is weakly closed, we have that $G(v)$ is weakly compact
and from Theorem \ref{kkm} it follows that $ \displaystyle \cap_{v
\in K} G(v) \neq \emptyset$, therefore from Remark \ref{remark2}
it follows that the problem (QHS) has a
solution.\\

\noindent{\bf Second Proof.} Using Remark \ref{remark2} we prove
the existence of an element $u \in K$ such that for every $v \in
K$ we have
$$A(Tu, Tv-Tu) + F(u, v-u)+ J^{\circ}(Su;Sv-Su) \geq 0.$$
In this case $u \in K$ will be the solution of systems (QHS).

We argue by contradiction. Let us assume that for each $u \in K$,
there exists $v \in K$ such that
\begin{equation}\label{masodik_ellen_egyenlet}
A(Tu, Tv-Tu) + F(u, v-u)+ J^{\circ}(Su;Sv-Su) < 0.
\end{equation}
 Now, we define the set-valued mapping $G: K \hookrightarrow
2^{K}$ by
\begin{equation}\label{masodik_egyenlet}
G(u)=\{v \in K \ : \  A(Tu, Tv-Tu) + F(u, v-u)+
J^{\circ}(Su;Sv-Su) < 0\}.
\end{equation}
From (\ref{masodik_ellen_egyenlet}) it follows that the set
$G(u)\neq \emptyset$ for every $u \in K$. Because the function
$A(\cdot, \cdot) + F(\cdot, \cdot)+ J^{\circ}(\cdot;\cdot)$ is
convex in the seconde variable, we get that $G(u)$ is a convex
set. Now, we prove that for every $v \in K$, the set $G^{-1}(v)=
\{u \in K \ : \ v \in G(u)\}$ is weakly open. Indeed, from weakly
upper semicontinuity of the function

$$u \mapsto A(Tu, Tv-Tu) + F(u, v-u)+
J^{\circ}(Su;Sv-Su)$$ it follows that
$$[G^{-1}(v)]^{c}= \{u \in K \ : \ A(Tu, Tv-Tu) + F(u, v-u)+
J^{\circ}(Su;Sv-Su) \geq 0\}$$ is weakly closed, therefore
$G^{-1}(v)$ is weakly open.\\
Now we verify iii) from Theorem \ref{Tarafdar_Fix}, i.e.
$\displaystyle \bigcup_{v \in K} G^{-1}(v)=K$. Because for every
$v \in K$ we have $G^{-1}(v) \subset K$, it follows that
$\displaystyle \bigcup_{v \in K} G^{-1}(v)\subset K$. Conversely,
let $u\in K$ be fixed. Since $G(u) \neq \emptyset$ there exists
$v_{0} \in K$ such that $v_{0} \in G(u)$. In the next step we
verify iv) Theorem \ref{Tarafdar_Fix}. We assert that $D=
\displaystyle \bigcap _{v \in K} [G^{-1}(v)]^{c}$ is empty or
weakly compact. Indeed, if $D\neq \emptyset$, then $D$ is a weakly
closed set of $K$ since it is the intersection of weakly closed
sets. But $K$ is weakly compact hence we get that $D$ is weakly
compact. Taking $O_{v}=G^{-1}(v)$ and $K_{0}=K_{1}=K$ we can apply
Theorem \ref{Tarafdar_Fix} to conclude that there exists $u_{0}
\in K$ such that $u_{0} \in G(u_{0})$. This give
$$0=A(Tu_{0}, Tu_{0}- Tu_{0}) + F(u_{0}, u_{0}-u_{0}) + J^{0}(Su_{0}; Su_{0} - Su_{0})<0, $$
which is a contradiction.Therefore the system (QHS) has a solution.\\

\begin{remark}\label{koerciv_remark}
If in Theorem \ref{Main_RES} the sets $K_{i}, \ i=\overline{1,n}$
are only convex and closed but not bounded we impose the following
coercivity condition.\\
{\bf (CC)} there exist $K_{i}^{0} \subset K_{i}$  compact sets and
$v_{i}^{0} \in K_{i}^{0}$ such that for all $v=(v_{1}, \ldots,
v_{n}) \in K_{1} \times \ldots \times K_{n} \setminus K_{1}^{0}
\times \ldots \times K_{n}^{0} $ we have
$$A(Tv, Tv^{0}-Tv) + F(v, v^{0}-v) + \displaystyle \sum_{i=1}^{n} J_{i}^{0}(Sv, S_{i}v_{i}^{0} - S_{i}v_{i}) <0,$$
where $v^{0}=(v_{1}^{0}, \ldots, v_{n}^{0})$. In this case the
problem {\bf (QHS)} has a solution.
\end{remark}

\section{Applications}

In this section we are concerned with two applications. In the
first application we study the relation between Nash equlibrum and
Nash generalized derivative equilibrum points for a
hemivariational inequalities system and in the second application
we give an existence result for an abstract class of
hemivariational inequalities systems.

Let $X_{1}, \ldots, X_{n}$ be Banach spaces and $K_{i} \subset
X_{i}$ and the functions $f_{i}: K_{1} \times \cdots \times K_{i}
\times \cdots \times K_{n} \rightarrow \mathbb{R}$ for $i\in \{1,
\ldots, n\}$. The following notion was introduced by J. Nash
\cite{N1}, \cite{N2}:

\begin{definition}\label{Nash_point}
An element $(u_{1}^{0}, \ldots, u_{n}^{0}) \in K_{1} \times \cdots
\times K_{n}$ is Nash equlibrum point of functions $f_{1}, \ldots,
f_{n}$ if for each $i \in \{1, \ldots, n\}$ and $(u_{1}, \ldots,
u_{n}) \in K_{1} \times \cdots K_{n}$ we have
$$f_{i}(u_{1}^{0}, \ldots, u_{i}, \ldots, u_{n}^{0}) \geq f_{i}(u_{1}^{0},
\ldots, u_{i}^{0}, \ldots, u_{n}^{0}).$$
\end{definition}

Now let  $D_{i} \subset X_{i}$ be open sets such that $K_{i}
\subset D_{i}$ for all $i \in \{1, \ldots, n\}$. We consider the
function $f_{i}: K_{1} \times \cdots \times D_{i} \times \cdots
K_{n} \rightarrow \mathbb{R}$ which are continuous and locally
Lipschitz in the i{\it th} variable. The next notion was
introduced recently by A. Krist\'{a}ly \cite{Kristaly3} and is a
little bit different form for functions defined on Riemannian
manifolds.

\begin{definition}
If $(u_{1}^{0}, \ldots, u_{n}^{0}) \in K_{1} \times \cdots \times
K_{n}$ is an element such that
$$f_{i}^{0}(u_{1}^{0}, \ldots, u_{n}^{0}; u_{i} - u_{i}^{0}) \geq 0,$$
for every $i=\{1, \ldots, n\}$ and $(u_{1}, \ldots, u_{n})\in
K_{1} \times \cdots \times K_{n}$ we say that $(u_{1}^{0}, \ldots,
u_{n}^{0})$ is a Nash generalized derivative points for the
functions $f_{1}, \ldots, f_{n}$.

\end{definition}

\begin{remark}
If the functions $f_{i}, i \in \{1, \ldots, n\}$ are
differentiable in the {\rm i}$th$ variable, then the above notion
coincides with the Nash stationary point introduced in {\rm
\cite{KKP}}.

\end{remark}

\begin{remark}\label{Nash_pont_nash-der}
Is is easy to observe that any Nash equlibrum point is a Nash
generalized derivative point.
\end{remark}

The following result is an existence result for Nash generalized
derivative points and is an infinite-dimensional version of a
result from the paper \cite{Kristaly3}. Therefore, if in Theorem
\ref{Main_RES} we choose $F_{i}=0, i\in \{1, \ldots, n\}$ and
$J=0$ we obtain the following result.

\begin{theorem} (i) Let $Y_{1}, Y_{2}, \ldots, Y_{n}$ and $X_{1}, X_{2}, \ldots,
X_{n},$ be a reflexive Banach spaces and $T_{i}: X_{i} \rightarrow
Y_{i}$  compact, linear operators. We consider the closed, convex,
bounded sets $K_{i} \subset X_{i}$ and the functions $A_{i}: Y_{1}
\times \cdots \times Y_{n} \rightarrow \mathbb{R}, i={1, \ldots,
n}$ which are locally Lipschitz in the {\rm i}th variable and
satisfies the condition {\rm (A)}. In these conditions, there
exists $(u_{1}^{0}, \ldots, u_{i}^{0}, \ldots, u_{n}^{0}) \in
K_{1} \times \cdots \times K_{i} \times \cdots \times K_{n}$ such
that for all $i \in \{1, \ldots, n\}$ and $(u_{1}, \ldots, u_{i},
\ldots, u_{n}) \in K_{1} \times \cdots \times K_{i} \times \cdots
\times K_{n}$ we have
$$A_{i}^{0}((T_{1}u_{1}^{0}, \ldots, T_{i}u_{i}^{0}, \ldots, T_{n}u_{n}^{0}); T_{i}u_{i} - T_{i}u_{i}^{0}) \geq 0,$$
i.e. $(u_{1}^{0}, \ldots, u_{i}^{0}, \ldots, u_{n}^{0})$ is a Nash
generalized derivative points for the function $A_{i}, i \in \{1,
\ldots, n\}$.\\

 (ii) If the sets $K_{i}, i=\{1, \ldots, n\}$ are only closed
and convex we suppose that there exists the bounded, closed sets
$K_{i}^{0} \subset K_{i}$ and $v_{i}^{0} \in K_{i}^{0}, i=\{1,
\ldots, n\}$ such that for every $(u_{1}, \ldots, u_{n}) \in K_{1}
\times \cdots \times K_{n} \setminus K_{1}^{0} \times \cdots
\times K_{n}^{0}$ we have
$$ A(Tu, Tv^{0}-Tu) < 0.$$
Then there exist $u^{0}=(u_{1}^{0}, \ldots, u_{i}^{0}, \ldots,
u_{n}^{0}) \in K_{1} \times \cdots \times K_{i} \times \cdots
\times K_{n}$ such that for all $i \in \{1, \ldots, n\}$ and
$u=(u_{1}, \ldots, u_{i}, \ldots, u_{n}) \in K_{1} \times \cdots
\times K_{i} \times \cdots \times K_{n}$ we have
$$A_{i}^{0}(Tu_{0}; T_{i}u_{i} - T_{i}u_{i}^{0}) \geq 0,$$
i.e. $u^{0}=(u_{1}^{0}, \ldots, u_{i}^{0}, \ldots, u_{n}^{0})$ is
a Nash generalized derivative points for the functions $A_{i}, i
\in \{1,
\ldots, n\}$.\\

\end{theorem}

In the next step we give an existence result for a general system
of hemivariational inequalities. In this case in Theorem
\ref{Main_RES} we choose $Y_{i}=Z_{i}, i \in \{1, \ldots, n\}$ and
we suppose that the functions $A_{i}: Y_{1} \times \cdots \times
Y_{i} \times \cdots \times Y_{n} \rightarrow \mathbb{R}$ are
differentiable in the i{\it th} variable for $i \in \{1, \ldots,
n\}$. In this case we suppose that the functions $A_{i}^{\prime}:
Y_{1} \times \cdots \times Y_{i} \times \cdots \times Y_{n} \times
Y_{i} \rightarrow \mathbb{R}$ are continuous for $i\in\{1, \ldots,
n\}$. Let also $J: Y_{1} \times \cdots \times Y_{i} \times \cdots
\times Y_{n} \rightarrow \mathbb{R}$ a locally Lipschitz regular
function.

Under these conditions we have the following result.

\begin{corollary}\label{hemivar_exist}
Let $J, A_{i}: Y_{1} \times \cdots \times Y_{i} \times \cdots
\times Y_{n} \rightarrow \mathbb{R}$ be the function as above and
suppose that the condition {\bf (TS)} holds and let $K_{i} \subset
X_{i}, i=\{1,\ldots, n\}$ be bounded, closed and convex sets.
Under these conditions there exist an element $u^{0}=(u_{1}^{0},
\ldots, u_{n}^{0}) \in K_{1} \times \cdots \times K_{n}$ such that
for every $u=(u_{1}, \ldots, u_{n}) \in K_{1} \times \cdots \times
K_{n}$ and $i \in\{1, \ldots, n\}$ we have:
$$A_{i}^{\prime}(Tu^{0}; T_{i}u_{i}- T_{i}u_{i}^{0}) + J_{i}^{0}(Tu^{0};T_{i}u_{i}- T_{i}u_{i}^{0} ) \geq 0.$$

\end{corollary}

If in Theorem \ref{Main_RES} we take $A_{i}=0$  then we obtain the
following existence result for a general class of hemivariational
inequalities systems.
\begin{corollary}\label{cor_exist_hem}
Let $K_{i} \subset X_{i}$ bounded, closed and convex  subsets of
the reflexive Banach spaces $X_{i}$ for $i\in \{1,\dots, n\}$. We
suppose that $F_{i}: K_{1} \times \cdots \times K_{n} \rightarrow
X_{i}^{\star}$ satisfies the condition {\bf (F)} and $J: Z_{1}
\times \cdots \times Z_{n} \rightarrow \mathbb{R}$ is a regular
locally Lipschitz function and the condition {\bf (TS)} holds.
Then there exists $u^{0}=(u_{1}^{0}, \ldots, u_{i}^{0}, \ldots,
u_{n}^{0}) \in K_{1} \times \cdots \times K_{i} \times \cdots
K_{n}$ such that for every $u=(u_{1}, \ldots, u_{i}, \ldots,
u_{n}) \in K_{1} \times \cdots \times K_{i} \times \cdots K_{n}$
and $i \in \{1, \ldots, n\}$ we have
$$\langle F_{i}(u^{0}); u_{i} - u_{i}^{0} \rangle _{i} +
J_{i}^{0}(Su^{0};S_{i}u_{i}- S_{i}u_{i}^{0} ) \geq 0.$$
\end{corollary}

The above result generalize the main result from the paper of A.
Krist\'{a}ly \cite{Kristaly_H}.

Indeed, let $\Omega \subset \mathbb{R}^{N}$ be a bounded, open
subset. Let $j: \Omega \times
\displaystyle\underbrace{\mathbb{R}^{k} \times \cdots
\mathbb{R}^{k}}_{n} \rightarrow \mathbb{R}$ a Carath\'{e}odory
function such that $j(x, \cdot, \ldots, \cdot)$ is locally
Lipschitz for every $x \in \Omega$ and satisfies the following
assumptions for all $i \in \{1, \ldots, n\}$: \\
{\rm $(j_{i})$} there exists $h_{1}^{i} \in
L^{\frac{p}{p-1}}(\Omega, \mathbb{R}_{+})$ and $h_{2}^{i} \in
L^{\infty}(\Omega, \mathbb{R}_{+})$ such that
$$|z_{i}| \leq h_{1}^{i}(x) + h_{2}^{i}(x)|y|^{p-1}_{\mathbb{R}^{kn}}$$
for almost $x \in \Omega$ and every $y=(y_{1}, \ldots, y_{n}) \in
\displaystyle \underbrace{\mathbb{R}^{k} \times \cdots \times
\mathbb{R}^{k}}_{n}$ and $z_{i} \in \partial_{i}j(x, y_{1},
\ldots, y_{n})$.

In this case let $S=(S_{1}, \ldots, S_{n}): X_{1} \times \cdots
\times X_{n} \rightarrow L^{p}(\Omega, \mathbb{R}^{k}) \times
\cdots L^{p}(\Omega, \mathbb{R}^{k})$ and $J\circ S: K_{1} \times
\cdots \times K_{n} \rightarrow \mathbb{R}$ is defined by
$$J(Su)= \displaystyle \int_{\Omega} j(x, S_{1}u_{1}(x), \ldots S_{n}u_{n}(x))dx.$$
Using a result from Clarke \cite{Cl} we have:
$$J_{i}^{0}(Su; S_{i}v_{i}) \leq \displaystyle \int_{\Omega} j_{i}^{0}(x, S_{1}u_{1}(x), \ldots S_{n}u_{x}; S_{i}v_{i}(x))dx, \leqno{({\rm I})}$$
for every $i\in \{1, \ldots, n\}$ and $v_{i} \in X_{i}$.\\
Therefore we have the following existence result obtained by
Krist\'{a}ly \cite{Kristaly_H}.

\begin{corollary}\label{Kr-hemi-tet}
Let $K_{i} \subset X_{i}$ bounded, closed and convex  subsets of
the reflexive Banach spaces $X_{i}$ for $i\in \{1,\dots, n\}$. We
suppose that $F_{i}: K_{1} \times \cdots \times K_{n} \rightarrow
X_{i}^{\star}$ satisfies the condition {\bf (F)} and $j: \Omega
\times \displaystyle\underbrace{\mathbb{R}^{k} \times \cdots
\mathbb{R}^{k}}_{n} \rightarrow \mathbb{R}$ a Carath\'{e}odory
function such that $j(x, \cdot, \ldots, \cdot)$  is a regular,
locally Lipschitz function satisfying condition $(j_{i})$ and the
condition {\bf (TS)} holds. Then there exists $u^{0}=(u_{1}^{0},
\ldots, u_{i}^{0}, \ldots, u_{n}^{0}) \in K_{1} \times \cdots
\times K_{i} \times \cdots K_{n}$ such that for every $u=(u_{1},
\ldots, u_{i}, \ldots, u_{n}) \in K_{1} \times \cdots \times K_{i}
\times \cdots K_{n}$ and $i \in \{1, \ldots, n\}$ we have
$$\langle F_{i}(u^{0}); u_{i} - u_{i}^{0} \rangle _{i} +
\displaystyle \int_{\Omega} j_{i}^{0}(x, S_{1}u_{1}^{0}(x),
\ldots, S_{n}u_{n}^{0}(x);S_{i}u_{i}(x)- S_{i}u_{i}^{0}(x))dx \geq
0.$$

\end{corollary}

\begin{remark}\label{FPRad_tetel}
If $n=1$ we obtain a similar result from the paper of
Panagiotopoulos, Fundo and R\u{a}dulescu {\rm \cite{PFRad_JOGO}}.
\end{remark}

\begin{remark}\label{nemkorlatos}
If the Banach spaces $X_{i}, i \in \{1, \ldots, n\}$ are separable
and the domain $\Omega \subset \mathbb{R}^{N}$ is unbounded then a
similar inequality to {\rm (I)} was proved in the paper D\'{a}lyai
and Varga \cite{DV}. Therefore, we can state a similar result as
Corollary \ref{Kr-hemi-tet} in the case when $\Omega \subset
\mathbb{R}^{N}$ is an unbounded domain.
\end{remark}
\textbf{Acknowledgement:} Cs. Varga has been supported by Grant
CNCSIS PN II ID PCE 2008 nr. 501, ID 2162. Both authors were
supported by Slovenian Research Agency grants No. P1-0292-0101 and
J1-2057-0101.


\begin{thebibliography}{99}

\bibitem{BBR} G. Bonanno,  G.M. Bisci, D. O'Regan,
{\it Infinitely many weak solutions for a class of quasilinear elliptic
systems}, Mathematical and Computer Modelling {\bf 52} (2010),
152-160.

\bibitem{BV} B.E. Breckner, Cs. Varga, 
{\it A multiplicity result for gradient-type systems with non-differentiable
term}, Acta Math. Hungar. {\bf 118} (2008), 85-104.

\bibitem{BHV} B.E. Breckner, A. Horv\'ath, Cs. Varga,
{\it A multiplicity result for a special class of
gradient-type systems with non-differentiable term}, Nonlinear
Anal. TMA, {\bf 70}(2009), 606-620.

\bibitem{KCC} K.-C. Chang, 
{\it Methods in Nonlinear Analysis,} Springer Verlag, Berlin, 2005.

\bibitem{Cl} {F. H. Clarke}, 
\textit{Optimization and Nonsmooth Analysis,}
SIAM, Philadelphia 1990.

\bibitem{DV} Z. D\'alyai, Cs. Varga,
{\it An existence result for hemivariational inequalities}, Electron. J. Differential
Equations {\bf 37} (2004), 1-17.

\bibitem{Fan_fix} K. Fan, 
\textit{Some properties of convex sets related to fixed point
theorem}, Math. Ann. {\bf 266} (1984), 519-537.

\bibitem{PFRad_JOGO} P.D. Panagiotopoulos, M. Fundo, V.
R\u{a}dulescu, {\it Existence Theorems of Hartman-Stampacchia Type
for Hemivariational Inequalities and Applications}, J.
 Global Optim. {\bf 15} (1999), 41-54.

\bibitem{KKP} G. Kassay, J. Kolumb\'{a}n, Zs. P\'{a}les,
{\it On Nash Stationary
Points}, Publ. Math. Debrecen {\bf 54} (1999), 267-279.

\bibitem{Kristaly_H}  A. Krist\' aly, {\it Hemivariational inequalities systems and
applications}, Mathematica {\bf 46} (2004), No.2, 161-168.


\bibitem{Kristaly1} A. Krist\' aly, {\it An existence result for gradient-type systems with a
nondifferentiable term on unbounded strips}, J. Math. Anal. Appl.
{\bf  229} (2004), 186-204.

\bibitem{Kristaly2} A. Krist\' aly, {\it  Existence of two non-trivial solutions
for a class of quasilinear elliptic variational systems on
stripe-like domains}, Proc. Edinburgh Math. Soc. {\bf 48} (2005),
1-13.

\bibitem{Kristaly3} A. Krist\' aly, {\it Location of Nash equilibria: a Riemannian
geometrical approach},  Proc. Amer. Math. Soc. {\bf 138} (2010),
1803-1810.

\bibitem{KVar2} A. Krist\' aly, Cs. Varga, {\it Variational-Hemivariational
 Inequalities on Unbounded Domains},
Studia Univ. Babe\c s-Bolyai, Mathematica, {\bf LV}(2010), Nr. 2,
3-87.


\bibitem{KRV} A. Krist\'aly, V. R\u adulescu, and Cs. Varga, {\it Variational
Principles in Mathematical Physics, Geometry, and Economics:
Qualitative Analysis of Nonlinear Equations and Unilateral
Problems}, Encyklopedia of Mathematics (No.~136), Cambridge
University Press, Cambridge, 2010.

\bibitem{LV} H. Lisei, Cs. Varga, {\it Multiple solutions for gradient elliptic systems with
nonsmooth boundary conditions}, Mediterr. J. Math. {\bf 8} (2011), 69-79.

\bibitem{MR} D. Motreanu, V. R\u adulescu,
\textit{ Variational  and non-variational methods in nonlinear
analysis and boundary value problems. Nonconvex Optimization and
its Applications}. Kluwer Academic Publishers, Dordrecht 2003.
\bibitem{MP} D. Motreanu, P. D. Panagiotopoulos,  {\it Minimax
Theorems and Qualitative Properties of the Solutions of
Hemivariational Inequalities,} Kluwer Academic Publishers,
Dordrecht 1999.

\bibitem{N1} J. Nash, {\it Equilibrium points in $n$-person games,}
{Proc. Nat. Acad. Sci. USA}, {\bf 36}(1950), 48-49.

\bibitem{N2} J.F. Nash, {\it Non-cooperative games.} Ann. of Math. (2) {\bf 54} (1951),
286-295.

\bibitem{PFR} P. D. Panagiotopoulos, M. Fundo,  V.
R\u adulescu, {\it Existence Theorems of Hartman-Stampacchia Type
for Hemivariational Inequalities and Applications}, Journal of
Global Optimization {\bf 15}(1999), 41-54.

\bibitem{Tar_Fix_eq} E. Tarafdar, \textit{Five equivalent theorems on a convex subset of a topological
vector space}, Commentationes Math. Univ. Carolinae, Vol. {\bf 30}
(1989), No.2, 323-326.

\bibitem{TarFix} E. Tarafdar, {\it A fixed point theorem equivalent to the Fan-Knaster-Kuratwski-Mazurkiewicz
theorem}, Journal Math. Anal. Appl. {\bf 128}(1987), 475-497.

\end{thebibliography}
\end{document}